\documentclass[12pt]{article}

\usepackage[english]{babel}
\usepackage{amsmath,amsthm,amscd,amssymb,latexsym,enumerate}
\usepackage[T1]{fontenc}
\usepackage{cite}
\usepackage{humanist}
\usepackage[letterpaper,hmargin=0.95in,vmargin=0.90in]{geometry}

\usepackage{a4}
\usepackage{amsmath}
\usepackage{amssymb}
\usepackage{amsthm}
\newtheorem{thm}{Theorem}[section]
\newtheorem{cor}[thm]{Corollary}
\newtheorem{lem}[thm]{Lemma}

\newtheorem{rem}[thm]{Remark}

\renewcommand{\b}{\beta}
\newcommand{\la}{\lambda}

\newcommand{\ind}{\operatorname{ind}}

\long\def\comment#1{{}}

\title{A question by Chihara about shell polynomials and indeterminate
moment problems\footnote{The first author was supported by grant No. 09-063996 and the second author by  Steno grant No.
09-064947 both from the Danish Research Council for Nature and Universe.}}

\author{Christian Berg\footnote{Corresponding author}, Jacob Stordal Christiansen}
\begin{document}
\maketitle

\begin{abstract}
The generalized Stieltjes--Wigert polynomials depending on parameters $0\le p<1$ and $0<q<1$ are discussed.
By removing the mass at zero of the N-extremal solution concentrated in the zeros of the $D$-function from
the Nevanlinna parametrization, we obtain a discrete measure $\mu^M$ which is uniquely determined by its
moments. We calculate the coefficients of the corresponding orthonormal polynomials $(P^M_n)$.
As noticed by Chihara, these polynomials are the shell polynomials corresponding to the maximal parameter
sequence for a certain chain sequence. We also find the minimal parameter sequence, as well as the parameter
sequence corresponding to the generalized Stieltjes--Wigert polynomials, and compute the value of related
continued fractions.
The mass points of $\mu^M$ have been studied in recent papers of Hayman, Ismail--Zhang and Huber.
In the special case of $p=q$, the maximal parameter sequence is constant and the determination of $\mu^M$
and $(P^M_n)$ gives an answer to a question posed by Chihara in 2001.
\end{abstract}


\section{Introduction}

In \cite{Ch2}, Chihara formulated an open problem concerning kernel
polynomials and chain sequences motivated by results in his old paper
\cite{Ch} and his monograph \cite{Ch1}. To formulate the problem
precisely, we need some notation and explanation, but roughly speaking
it deals with the following observation of Chihara.

Let $(k_n)$ denote the kernel polynomials of an indeterminate
Stieltjes moment problem. The corresponding shell polynomials
$(p_n^h)$, parametrized by the initial condition $0<h_0\le M_0$ for
the non-minimal parameter sequences $h=(h_n)$ of the associated  chain
sequence, are orthogonal with respect to the measure
$$
\mu^h = \mu^M + (M_0/h_0-1)\mu^M(\mathbb R)\delta_0.
$$
In the case of the generalized Stieltjes--Wigert polynomials
$S_n(x;p,q)$ with $p=q$, Chihara observed that the maximal parameter
sequence is constant
\[
M_n=\frac{1}{1+q}
\]
and for this special case Chihara's question is:

\lq\lq {\it Find the measure $\mu^M$ which has the property that the
  Hamburger moment problem is determinate, but if mass is added at the
  origin, the Stieltjes problem becomes indeterminate}\rq\rq

In this paper we find the measure $\mu^M$ as the discrete measure
$$
\mu^M=\sum_{n=1}^\infty \rho_n\delta_{\tau_n}
$$
obtained by removing the mass at zero from an N-extremal solution to
the generalized Stieltjes--Wigert moment problem, and the numbers
$\tau_n$ behave like
\[
\tau_n=q^{-2n-1/2}\bigl(1+\mathcal{O}(q^n)\bigr) \mbox{ as }\,  n\to\infty.
\]
For $p$, $q$ small enough or $n$ sufficiently large, there are constants $b_j$, $j\geq 1$, such that $\tau_n$ is given by
$$
\tau_n=q^{-2n-1/2}\Bigl(\sum_{j=1}^\infty b_j q^{jn}\Bigr),
$$
see Theorem~\ref{prop:zeros} for details.
These results are due to Hayman \cite{Hay}, Ismail--Zhang \cite{I:Z}, and Huber \cite{Hub}.
It does not seem possible to find more explicit formulas for the numbers
$\tau_n$ because this is equivalent to finding the zeros of the
$q$-Bessel function $J_\nu^{(2)}(z;q)$.

We also find explicit formulas for the coefficients of the orthonormal
polynomials associated with the measure $\mu^M$, see Theorem~\ref{thm:GSW},
and compute the minimal and maximal parameter sequences as well as the parameter
sequence corresponding to $S_n(x;p,q)$ in Theorem~\ref{thm:chaingenSW}.
The explicit expressions at hand allow us to show that
\[
1-\cfrac{\beta_{n}}{1-\cfrac{\beta_{n+1}}{1-\dotsb}}
=\frac{q\bigl((pq^{n-1};q)_\infty-(q^{n-1};q)_\infty \bigr)}
{\bigl(1+q-(1+p)q^{n}\bigr)
\bigl((pq^{n};q)_\infty-(q^{n};q)_\infty\bigr)}
\]
for every $n\geq 1$, where
\[
\beta_n=\frac{q(1-q^n)(1-pq^n)}
{\bigl(1+q-(1+p)q^n\bigr)\bigl(1+q-(1+p)q^{n+1}\bigr)}.
\]

\section{Preliminaries}
It is well known that chain sequences can be used to characterize
those three-term recurrence relations for orthogonal polynomials which
has a measure of orthogonality supported by $[0,\infty[$, cf. \cite {Ch1}.
The moments of such a measure is called a Stieltjes moment sequence,
and it is called determinate in the sense of Stieltjes (in short
det(S)) if there is only one measure supported on $[0,\infty[$ with
these moments, while it is called indeterminate in the sense of
Stieltjes (in short indet(S)) if there are different measures on the
half-line with these moments.

If a Stieltjes moment sequence is indet(S), then there are also
measures with the same moments and not
supported by the positive half-line. This follows from the Nevanlinna
parametrization of the indeterminate Hamburger moment problem. If the
Stieltjes moment sequence is det(S), it is still possible that it is
an indeterminate Hamburger moment sequence. See \cite{B:V} for
concrete examples.

In the following, let $(p_n)$ be a sequence of
monic orthogonal polynomials for a positive measure $\mu$ with
moments of any order and infinite support contained in $[0,\infty[$.
We  denote by $(k_n)$ the sequence of monic orthogonal polynomials
with respect to the measure $xd\mu(x)$. They are called \emph{kernel
polynomials} because they are the monic version of the reproducing
kernels
$$
K_n(x,y)=\sum_{k=0}^n p_k(x)p_k(y)/\Vert p_k \Vert^2,\quad \Vert p_k \Vert^2=\int
p^2_k(x)\,d\mu(x)
$$
when $y=0$, i.e.,
$$
k_n(x)=\frac{\Vert p_n \Vert^2}{p_n(0)}K_n(x,0).
$$
The three-term recurrence relation for the kernel polynomials is given
as
\begin{equation}
\label{eq:3t}
k_n(x)=(x-d_n)k_{n-1}(x)-\nu_nk_{n-2},\quad n\ge 1
\end{equation}
(with the convention that $k_{-1}=0$, $\nu_1$ is not defined). It is known,
cf. \cite{Ch}, that
\begin{equation}\label{eq:chainseq}
\beta_n=\nu_{n+1}/(d_nd_{n+1}),\quad n\ge 1
\end{equation}
is a chain sequence which
does not determine the parameter sequence uniquely. In this case there exists
a largest
$M_0>0$ such that for any $0\le h_0\le M_0$, there is a parameter
sequence $h_n$, $n\ge 0$, such that
\begin{equation}\label{eq:par}
\beta_n=h_n(1-h_{n-1}),\quad n\ge 1.
\end{equation}
The parameter sequence $M_n=h_n$ (resp. $m_n=h_n$) determined by
$h_0=M_0$ (resp. $h_0=m_0=0$) is called the maximal (resp. minimal)
parameter sequence. For each parameter sequence $h=(h_n)$
with $0<h_0\le M_0$, there exists a family of monic orthogonal polynomials
$(p_n^h)$ on $[0,\infty[$ which all have $(k_n)$ as kernel
polynomials. The polynomials $(p_n^h)$ are called the {\it shell
  polynomials} of  the kernel polynomials $(k_n)$. The coefficients in the
three-term recurrence relation
\begin{equation}\label{eq:pnh}
p_n^h(x)=(x-c_n^h)p_{n-1}^h(x)-\lambda_n^hp_{n-2}^h(x)
\end{equation}
are given explicitly in \cite{Ch} in terms of $d_n,h_n$ by
\begin{equation}
\label{eq:pnh1}
c_1^h=h_0d_1, \quad
c_{n+1}^h=(1-h_{n-1})d_n+h_nd_{n+1}, \quad n\ge 1
\end{equation}
and
\begin{equation}
\label{eq:pnh2}
\lambda_{n+1}^h=(1-h_{n-1})h_{n-1}d_n^2, \quad n\ge 1.
\end{equation}

Theorem 2 in \cite{Ch} states:
\begin{thm}
\label{thm:chain} The polynomials $(p_n^M)$ are orthogonal with
respect to a determinate measure $\mu^M$
 which has no mass at 0.

For $0<h_0<M_0$, the polynomials $(p_n^h)$ are orthogonal with respect
to
\begin{equation}\label{eq:chain1}
\mu^h=\mu^M+(M_0/h_0-1)\mu^M(\mathbb R)\delta_0,
\end{equation}
where $\delta_0$ denotes the Dirac measure with mass 1 at 0.

The measure $\mu^h$ is indet(S) if and only if $xd\mu^h(x)=xd\mu^M(x)$
is indet(S).
\end{thm}

\begin{rem}
{\rm Recall that for a measure $\mu$, the proportional measure
  $\lambda\mu$ ($\lambda>0$) leads to the same
 monic orthogonal polynomials as $\mu$. The normalization in \eqref{eq:chain1} is chosen so that $\lambda\mu^h$ precisely corresponds to $\lambda\mu^M$ for any $\lambda>0$.}
\end{rem}

In all of this paper we shall be focusing on the case where
$xd\mu^M(x)$ is indet(S), i.e., when the kernel polynomials correspond
to an indeterminate Stieltjes moment problem.

Concerning the ``if and only if'' statement of the theorem,
it is easy to see that if $\mu^h$ is indet(S), then $xd\mu^h(x)$ is
indet(S). The reverse implication is proved in \cite[p.~6--7]{Ch}. It
is also a consequence of \cite[Lemma 5.4]{B:T}.

The measure $\mu^M$ is determinate in the sense of Hamburger and
$xd\mu^M(x)$ is indet(S). Using the terminology of \cite[Sect.~5]{B:T}, we see
that the index of determinacy $\ind(\mu^M)$ is $0$. The
measures on $[0,\infty[$ of index zero were characterized in
\cite[Thm.~5.5]{B:T} as the discrete measures $\sigma$
defined in the following way: Take any Stieltjes moment sequence $(s_n)$
which is indet(S) and let $\nu_0$ be the corresponding N-extremal solution
which has a mass at 0. Define $\sigma$ by
\[
\sigma=\nu_0-\nu_0(\{0\})\delta_0.
\]
In other words, if $(P_n)$ are the orthonormal polynomials corresponding
to $(s_n)$ and if
\begin{equation}\label{eq:D}
D(z)=z\sum_{n=0}^\infty P_n(z)P_n(0),
\end{equation}
then $D$ has simple zeros $\tau_0=0<\tau_1<\ldots<\tau_n<\ldots$ and
\begin{equation}\label{eq:sigma}
\nu_0=\sum_{n=0}^\infty \rho_n\delta_{\tau_n},\quad
\sigma=\sum_{n=1}^\infty \rho_n\delta_{\tau_n},
\end{equation}
where
\begin{equation}\label{eq:rho1}
\rho_n^{-1}=\sum_{k=0}^\infty P_k^2(\tau_n).
\end{equation}

Already Stieltjes observed that removing the mass at zero of the
solution $\nu_0$ to an indeterminate Stieltjes problem leads to a
determinate solution, see \cite[Sect.~65]{St}. This phenomenon was
exploited in \cite{B:C} for indeterminate Hamburger moment problems
and carried on in Berg--Dur{\'a}n \cite{B:D}. It follows that all the
measures $\mu^h$ given by \eqref{eq:chain1} for $0<h_0<M_0$ are N-extremal.

\section{The generalized Stieltjes--Wigert polynomials}

For $0<q<1$ and $0\le p<1$, we consider the moment sequence
\begin{equation}
s_n=(p;q)_nq^{-(n+1)^2/2},\quad n\ge 0
\end{equation}
given by the integrals
\begin{equation}\label{eq:SW1}
\frac{1}{\sqrt{2\pi\log(1/q)}}\int_0^\infty x^n\exp\left(-\frac{(\log
x)^2}{2\log(1/q)}\right)(p,-p/\sqrt{q}x;q)_\infty\,dx.
\end{equation}
We call it the generalized Stieltjes--Wigert moment sequence because
it is associated with the generalized Stieltjes--Wigert polynomials
\begin{equation}\label{eq:genSW1}
S_n(x;p,q)=(-1)^nq^{-n(n+1/2)}(p;q)_n\sum_{k=0}^n
\left[\begin{matrix}n\\k\end{matrix}\right]_q\frac{q^{k^2}(-\sqrt{q}x)^k}
{(p;q)_k},
\end{equation}
where we follow the monic notation and normalization of \cite[p.~174]{Ch1}
for these polynomials.
We have used the Gaussian $q$-binomial coefficients
$$
\left[\begin{matrix}n\\k\end{matrix}\right]_q=
\frac{(q;q)_n}{(q;q)_k(q;q)_{n-k}},
$$
involving the $q$-shifted factorial
$$
(z;q)_n=\prod_{k=1}^n(1-zq^{k-1}),\quad z\in\mathbb C, \;\, n=0,1,\ldots,\infty.
$$
We refer to \cite{G:R} for information about this notation and
$q$-series.

The Stieltjes--Wigert polynomials corresponds to the special case $p=0$.
In his famous memoir \cite{St}, Stieltjes noticed that the special
values $\log(1/q)=1/2$ and $p=0$ give an example of an indeterminate moment
problem, and Wigert \cite{Wig} found the corresponding orthonormal
polynomials.  The normalization is the same as in Szeg\H{o}
\cite{Sz1}. Note that
\begin{equation}
s_0=1/\sqrt{q}.
\end{equation}
The Stieltjes--Wigert moment problem has been extensively
studied in \cite{Chr} using a slightly different normalization.

For the generalized Stieltjes--Wigert polynomials, the orthonormal
version is given as
\begin{equation}
\label{eq:SW2}
P_n(x;p,q)=(-1)^nq^{{n}/{2}+{1}/{4}}\sqrt{\frac{(p;q)_n}{(q;q)_n}}\sum_{k=0}^n
(-1)^k
\left[\begin{matrix}n\\k\end{matrix}\right]_q\frac{q^{k^2+{k}/{2}}}{(p;q)_k}x^k.
\end{equation}
From \eqref{eq:SW2} we get
\begin{equation}
\label{eq:SW3}
P_n(0;p,q)=(-1)^nq^{n/2+1/4}\sqrt{\frac{(p;q)_n}{(q;q)_n}}
\end{equation}
and hence, by the $q$-binomial theorem,
\begin{equation}\label{eq:SW4}
\sum_{n=0}^\infty P_n^2(0;p,q)=\sqrt{q}\sum_{n=0}^\infty
\frac{(p;q)_n}{(q;q)_n}q^{n}=\sqrt{q}\frac{(pq;q)_\infty}
{(q;q)_\infty}.
\end{equation}

From the general theory in \cite{Ak} we know that the generalized Stieltjes--Wigert moment
sequence has an N-extremal solution $\nu_0$ which has the mass
\begin{equation}\label{eq:massat0}
c=\frac{(q;q)_\infty}{\sqrt{q}(pq;q)_\infty}
\end{equation}
($=$ the reciprocal of the value in \eqref{eq:SW4}) at 0.
It is a discrete measure concentrated at the zeros of the entire function
\begin{equation}
\label{eq:Dfct}
D(z)=z\sum_{n=0}^\infty P_n(0;p,q)P_n(z;p,q).
\end{equation}
The measure
$\tilde\mu=\nu_0-c\delta_0$ is determinate, cf., e.g., \cite[Thm.~7]{B:C}. The moment sequence $(\tilde s_n)$ of
$\tilde\mu$ equals the Stieltjes--Wigert moment sequence except for
$n=0$, 
$$
\tilde s_n=\left\{\begin{array}{ll}
q^{-1/2}\left[1-(q;q)_\infty/(pq;q)_\infty\right] & \mbox{if $n=0$},\\
(p;q)_nq^{-(n+1)^2/2} & \mbox{if $n\ge 1$},
\end{array}\right.
$$
and similarly the corresponding Hankel matrices $\mathcal
H$ and $\tilde{\mathcal H}$ only differ at the
entry $(0,0)$. The orthonormal polynomials associated with $(\tilde
s_n)$ will be denoted $\tilde P_n(x;p,q)$. We call them the {\it modified
generalized Stieltjes--Wigert polynomials} and they will be determined
in Section 4.

With \eqref{eq:SW2}--\eqref{eq:SW3} at hand, we can find the entire function $D$ in \eqref{eq:Dfct} explicitly.
The following generating function leads to the power series expansion of $D$.
\begin{lem}
For $\vert t\vert<1$, we have
\begin{equation}
\label{eq:genfct}
\sum_{n=0}^\infty\frac{(p;q)_n}{(q;q)_n}
\biggl(\sum_{k=0}^n\left[\begin{matrix}n\\k\end{matrix}\right]_q
\frac{q^{k^2+{k}/{2}}}{(p;q)_k}z^k\biggr) t^n=
\frac{(pt;q)_\infty}{(t;q)_\infty}
\sum_{n=0}^\infty\frac{q^{n^2+{n}/{2}}}{(pt,q;q)_n}(zt)^n.
\end{equation}
\end{lem}
\begin{proof}
Since the double series on the left-hand side is absolutely convergent, we can interchange the order of summation to get
\[
LHS=\sum_{k=0}^\infty\frac{q^{k^2+{k}/{2}}}{(p,q;q)_k}z^k
\sum_{n=k}^\infty\frac{(p;q)_n}{(q;q)_{n-k}}t^n.
\]
Shifting the index of summation on the inner sum, the $q$-binomial theorem leads to
\[
LHS=\frac{1}{(t;q)_\infty}\sum_{k=0}^\infty
\frac{q^{k^2+{k}/{2}}}{(q;q)_k}(ptq^k;q)_\infty(zt)^k.
\]
We thus arrive at \eqref{eq:genfct}.
\end{proof}
Set $t=q$ and replace $z$ by $-z$ in \eqref{eq:genfct} to get
\begin{equation}
\label{eq:explD}
D(z)=z\sqrt{q}\,\frac{(pq;q)_\infty}{(q;q)_\infty}
\sum_{n=0}^\infty(-1)^n\frac{q^{n(n+1)}}{(pq,q;q)_n}(z\sqrt{q})^n.
\end{equation}
This is essentially the $q$-Bessel function $J_\nu^{(2)}(z;q)$ for $q^\nu=p$, cf. \cite{G:R}.

Besides $\tau_0=0$, the zeros $\tau_n$ of \eqref{eq:explD} cannot be found explicitly.
However, the asymptotic behaviour of $\tau_n$ for $n$ large can be described up to a small error.
General results of Bergweiler--Hayman \cite{B:W} show that
\begin{equation}
\label{eq:zeros1}
\tau_n=Aq^{-2n}\bigl(1+\mathcal{O}(q^n)\bigr) \,\mbox{ as }\, n\to\infty
\end{equation}
for some constant $A>0$. In fact, $A=q^{-1/2}$ as follows from later work of Hayman. He proved in \cite{Hay} that
\begin{thm}
Given $k\geq 1$, there are constants $b_1,\ldots, b_k$ (depending on $p$, $q$) such that
\begin{equation}
\label{eq:zeros2}
\tau_n=q^{-2n-1/2}\Bigl(1+\sum_{j=1}^k b_jq^{jn}+\mathcal{O}\bigl(q^{(k+1)n}\bigr)\Bigr)
\,\mbox{ as }\, n\to\infty.
\end{equation}
\end{thm}
The first few values of the constants are
\begin{eqnarray*}
&&b_1=-\frac{1+p}{(1-q)\psi^2(q)}, \quad b_2=0, \\
&&b_3=-\frac{q(1+q^2)(1+p^3)+2pq(1+p)(1+q+q^2)}{(1-q)(1-q^2)(1-q^3)\psi^2(q)} \\
&&\qquad\,+ \frac{(1+p)^3}{(1-q)^3\psi^6(q)}\sum_{j=1}^\infty\frac{(2j-1)q^{2j-1}}{1-q^{2j-1}}, \\
&&b_4=b_1b_3,
\end{eqnarray*}
where
\[
\psi(q)=\sum_{n=0}^\infty q^{n(n+1)/2}=\frac{(q^2;q^2)_\infty}{(q;q^2)_\infty}.
\]
Even stronger results were recently obtained by Ismail--Zhang \cite{I:Z} and Huber \cite{Hub}.
They showed that for $n$ sufficiently large (in \cite{I:Z}) or for every $n$ when $p$, $q$ are small enough (in \cite{Hub}),
\begin{thm}
\label{prop:zeros}
There are constants $b_j$, $j\geq 1$, such that $\tau_n$ is given exactly by the convergent series
\begin{equation}
\label{eq:zeros3}
\tau_n=q^{-2n-1/2}\Bigl(1+\sum_{j=1}^\infty b_jq^{jn}\Bigr).
\end{equation}
\end{thm}
The $b_j$'s satisfy a somewhat complicated recursion formula that in principle allows for determining $b_{j+1}$ from $b_1,\ldots, b_j$. See \cite{Hub} for details.

\section{The modified generalized Stieltjes--Wigert polynomials}
It is a classical fact, cf. \cite[p.~3]{Ak}, that the
orthonormal polynomials $(P_n)$ corresponding to a moment sequence
$(s_n)$ are given by the formula
\begin{equation}\label{eq:detpol}
P_n(x)=\frac{1}{\sqrt{D_{n-1}D_n}}\det\begin{pmatrix} s_0 & s_1
  & \cdots & s_n\\ \vdots & \vdots & \ddots & \vdots\\
  s_{n-1}& s_n & \cdots & s_{2n-1}\\
1 & x & \cdots & x^n
\end{pmatrix},
\end{equation}
where
\[
D_n=\det \mathcal H_n, \quad \mathcal H_n=(s_{i+j})_{0\leq i,j \leq n}.
\]
In this way Wigert calculated the polynomials $P_n(x;0,q)$ and we shall
follow the same procedure for $P_n(x;p,q)$ and $\tilde P_n(x;p,q)$.
The calculation of $\tilde P_n(x;0,q)$ was carried out in \cite{B:S}.

It will be convenient to use the notation
\begin{equation}
\label{eq:rho}
\Delta_n:=(pq^n;q)_\infty-(q^n;q)_\infty, \quad n\geq 0.
\end{equation}
Writing
\begin{equation}\label{eq:swtilde}
P_n(x;p,q)=\sum_{k=0}^n b_{k,n}x^k,\quad \tilde P_n(x;p,q)=\sum_{k=0}^n \tilde
b_{k,n}x^k,
\end{equation}
we have:
\begin{thm}
\label{thm:GSW}
For $0\le k\le n$,
\begin{equation}\label{eq:swtilde1}
\tilde b_{k,n}=\tilde C_n (-1)^k\left[\begin{matrix}n\\k\end{matrix}\right]_q
\frac{q^{k^2+{k}/{2}}}{(p;q)_k}
\left[1-\frac{1-q^k}{1-pq^k}\frac{(q^{n+1};q)_\infty}{(pq^{n+1};q)_\infty}\right],
\end{equation}
where
\begin{equation}
\begin{aligned}
\label{eq:swtilde2}
\tilde C_n&=(-1)^nq^{{n}/{2}+{1}/{4}}\sqrt{\frac{(p;q)_n}{(q;q)_n}}
\left[\left(1-\frac{(q^n;q)_\infty}{(pq^n;q)_\infty}\right)\left(1-\frac{(q^{n+1};q)_\infty}{(pq^{n+1};q)_\infty}\right)\right]^{-1/2} \\
&=(-1)^nq^{{n}/{2}+{1}/{4}}\sqrt{\frac{(p;q)_{n+1}}{(q;q)_n}}\,
\frac{(pq^{n+1};q)_\infty}{\sqrt{\Delta_n\Delta_{n+1}}},
\end{aligned}
\end{equation}
i.e.,
\begin{equation}
\label{eq:swtilde3}
\tilde b_{k,n}=b_{k,n}
\left[(pq^{n+1};q)_\infty-\frac{1-q^k}{1-pq^k}(q^{n+1};q)_\infty\right]
\sqrt{\frac{1-pq^n}{\Delta_n\Delta_{n+1}}}
\end{equation}
Moreover,
\begin{equation}
\label{eq:swtilde4}
\tilde D_n=\frac{\Delta_{n+1}}{(pq^{n+1};q)_\infty}D_n,
\end{equation}
where $D_n=\det\mathcal H_n$ and $\tilde D_n=\det\tilde{\mathcal H}_n$.
\end{thm}

\begin{proof}
We first recall the Vandermonde determinant
\begin{equation}\label{eq:van}
V_n(x_1,\ldots,x_n)=\det\begin{pmatrix}
1 & 1 & \cdots & 1 \\
x_1 & x_2 & \cdots & x_n\\
\vdots & \vdots& \ddots& \vdots\\
x_1^{n-1} & x_2^{n-1} & \cdots & x_n^{n-1}
\end{pmatrix}
=\prod_{1\le i<j\le n}(x_j-x_i).
\end{equation}
Using the moments $s_n=(p;q)_nq^{-(n+1)^2/2}$, we find
\begin{equation}\label{eq:Dn1}
D_n=\left(\prod_{j=0}^n s_j\right)\det\left(s_{i+j}/s_j\right)=
\left(\prod_{j=1}^n(p;q)_j\right)q^{-\tfrac12\sigma_{n+1}}\det\left(s_{i+j}/s_j\right),
\end{equation}
where $\sigma_n=\sum_{j=0}^n j^2=n(n+1)(2n+1)/6$.
Noting that
\begin{equation}\label{eq:h1}
s_{i+j}/s_j=(pq^j;q)_i\, q^{-i^2/2}q^{-i(j+1)},
\end{equation}
we get
\begin{equation}
\label{eq:Dn2}
D_n=\left(\prod_{j=1}^n(p;q)_j\right)q^{-\tfrac12(\sigma_{n+1}+\sigma_n)}\det\left((pq^j;q)_i\,q^{-i(j+1)}\right).
\end{equation}

The last determinant can be simplified in the following way: Multiply
the first row (corresponding to $i=0$) by $p/q$ and add it to the
second row ($i=1$). Then the second row becomes $q^{-(j+1)},j=0,1,\ldots,n$, and the
determinant is not changed. The third row ($i=2$) has
the entries
$$
q^{-2(j+1)}-p(1+1/q)q^{-(j+1)}+{p^2}/{q}, \quad j=0,1,\ldots,n,
$$
so adding the first row multiplied by $-p^2/q$ and the second row
multiplied by $p(1+1/q)$ to the third row, changes the third row to
$q^{-2(j+1)}$, $j=0,1,\ldots,n$, and the determinant is not changed. If we go on like
this, we finally get
\begin{equation}\label{eq:Dn3}
D_n=\left(\prod_{j=1}^n(p;q)_j\right)q^{-\tfrac12(\sigma_{n+1}+\sigma_n)}\det\left(q^{-i(j+1)}\right).
\end{equation}
The last determinant is precisely $V_{n+1}(q^{-1},\ldots,q^{-(n+1)})$ and by \eqref{eq:van} equal to
\[
\prod_{i=1}^n\prod_{j=i+1}^{n+1}\left(q^{-j}-q^{-i}\right)=
\prod_{i=1}^n q^{-(n+1-i)(n+2+i)/2}(q;q)_{n+1-i}.
\]
After some reduction, we get
\begin{equation}
\label{eq:Dn4}
V_{n+1}(q^{-1},\ldots,q^{-(n+1)})=q^{-n(n+1)(n+2)/3}\prod_{j=1}^n(q;q)_j.
\end{equation}
Hence,
\begin{equation}\label{eq:Dn5}
D_n=\left(\prod_{j=1}^n(p,q;q)_j\right)q^{-(n+1)(2n+1)(2n+3)/6}
\end{equation}
and for later use we note that
\begin{equation}\label{eq:Dn/D}
D_n/D_{n-1}=(p,q;q)_n\,q^{-(2n+1)^2/2}.
\end{equation}

We denote by $A_{r,s}$ (resp.~$\tilde{A}_{r,s}$) the cofactor of entry
$(r,s)$ of the Hankel matrix $\mathcal H_n$ (resp.~$\tilde{\mathcal H}_n$), where $r,s=0,1,\ldots,n$.
(Note that entry $(r,s)$ is in row number $r+1$ and column number $s+1$.)
When $r=0$ or $s=0$, we clearly have $A_{r,s}=\tilde{A}_{r,s}$.
For $0\le s\le n$, we get
\begin{equation*}
\begin{aligned}
A_{n,s}=&(-1)^{n-s}\det\Bigl(s_{i+j}\;\Big\vert\;
{}^{\,\,i=0,\ldots,n-1}_{j=0,\ldots,n;\,j\neq s}\Bigr) \\
=&(-1)^{n-s}\biggl(\prod_{{}^{j=0}_{j\neq s}}^{n}s_j \biggr)
\det\Bigl({s_{i+j}}/{s_j}\;\Big\vert\; {}^{\,\,i=0,\ldots,n-1}_{j=0,\ldots,n; \,j\neq s}\Bigr) \\
=&(-1)^{n-s}\biggl(\prod_{{}^{j=0}_{j\neq s}}^{n}(p;q)_j\biggr)
q^{-\tfrac12\bigl(\sigma_{n+1}+\sigma_{n-1}-(s+1)^2\bigr)} \\
&\times
\det\Bigl((pq^j;q)_i\,q^{-i(j+1)}\;\Big\vert\; {}^{\,\,i=0,\ldots,n-1}_{j=0,\ldots,n; \,j\neq s}\Bigr).
\end{aligned}
\end{equation*}
However, the last determinant can be simplified like the
simplifications from \eqref{eq:Dn2} to \eqref{eq:Dn3} to give the
Vandermonde determinant $V_n(q^{-(j+1)}\;|\;j=0,\ldots,n, j\neq s)$.
To calculate this determinant, we observe that
\begin{equation*}
\begin{aligned}
V&_{n+1}(q^{-1},\ldots,q^{-(n+1)}) \\
  &=V_n(q^{-(j+1)}\;|\;j=0,\ldots,n,j\neq s)\prod_{j=0}^{s-1}(q^{-(s+1)}-q^{-(j+1)})
\prod_{j=s+1}^{n}(q^{-(j+1)}-q^{-(s+1)})\\
  &=V_n(q^{-(j+1)}\;|\;j=0,\ldots,n,j\neq
s)(q;q)_sq^{-s(s+1)}(q;q)_{n-s}q^{-\tfrac12(n-s)(n+s+3)}
\end{aligned}
\end{equation*}
and hence
\begin{equation}
\begin{aligned}\label{eq:Ans}
A&_{n,s} \nonumber \\
 &=\frac{(-1)^{n-s}}{(q;q)_n(p;q)_s}\prod_{j=0}^n(p;q)_j\left[\begin{matrix}n\\s\end{matrix}\right]_q
V_{n+1}(q^{-1},\ldots,q^{-(n+1)})q^{-\tfrac12(\sigma_{n+1}+\sigma_{n-1}-n(n+3)-1)}q^{s^2+s/2}\nonumber\\
 &=\frac{(-1)^{n-s}}{(q;q)_n(p;q)_s}\left[\begin{matrix}n\\s\end{matrix}\right]_q
D_n q^{(n+1)(n+1/2)}q^{s^2+s/2}.
\end{aligned}
\end{equation}
Using \eqref{eq:detpol} it is now easy to verify formula \eqref{eq:SW2} for
the generalized Stieltjes--Wigert polynomials $P_n(x;p,q)$.

Expanding after the first column, we get
$$
\tilde D_n=D_n-cA_{0,0},\quad c=\frac{(q;q)_\infty}{\sqrt{q}(pq;q)_\infty}
$$
and a calculation as above leads to
\begin{equation*}
\begin{aligned}
A_{0,0}&=\det\left(s_{i+j}\;|\;i,j=1,\ldots,n\right)  \\
 &=\left(\prod_{j=1}^n s_{j+1}\right)\det\left({s_{i+j}}/{s_{j+1}}\;|\;i,j=1,\ldots,n\right)\\
 &=\left(\prod_{j=2}^{n+1}(p;q)_j\right)q^{-\tfrac12(\sigma_{n+2}+\sigma_{n-1}-5)}V_n(q^{-j},j=3,\ldots,n+2)\\
 &=\left(\prod_{j=2}^{n+1}(p;q)_j\right)q^{-\tfrac12(\sigma_{n+2}+\sigma_{n-1}-5)}q^{-n(n-1)}V_n(q^{-j},j=1,\ldots,n).
\end{aligned}
\end{equation*}
Using \eqref{eq:Dn4} with $n$ replaced by $n-1$ and \eqref{eq:Dn5},
we find
$$
A_{0,0}=D_n\frac{(p;q)_{n+1}\sqrt{q}}{(1-p)(q;q)_n}
$$
which gives \eqref{eq:swtilde4}.

For $1\le s\le n$, we find
$$
\tilde A_{n,s}=A_{n,s}-c(-1)^{n-s}
\det\Bigl(s_{i+j}\;\Big\vert\; {}^{\,\,i=1,\ldots,n-1}_{j=1,\ldots,n;\,j\neq s}\Bigr)
$$
and the determinant on the right-hand side can be calculated by the same method as above
to be
\begin{equation*}
\begin{aligned}
\biggl(&\prod_{{}^{j=1}_{j\neq s}}^{n}s_{j+1}\biggr)
\det\Bigl({s_{i+j}}/{s_{j+1}}\;\Big\vert\; {}^{\,\, i=1,\ldots,n-1}_{j=1,\ldots,n;\,j\neq s}\Bigr) \\
&=\biggl(\prod_{{}^{j=1}_{j\neq s}}^{n}(p;q)_{j+1}\biggr)
q^{-\tfrac12(\sigma_{n+2}+\sigma_{n-2}-5-(s+2)^2)}
V_{n-1}(q^{-(j+2)},j=1,\ldots,n;j\neq s)\\
&=D_{n-1}\frac{(p;q)_n(p;q)_{n+1}}{(1-p)(p;q)_{s+1}(q;q)_{n-s}(q;q)_{s-1}}q^{-n^2-(n-1)/2+s(s+1/2)}.
\end{aligned}
\end{equation*}
This leads to
\begin{equation}\label{eq:Dn6}
\tilde
A_{n,s}=A_{n,s}\left[1-\frac{1-q^s}{1-pq^s}\frac{(q^{n+1};q)_\infty}{(pq^{n+1};q)_\infty}\right]
\end{equation}
which also holds for $s=0$ because then $\tilde A_{n,0}=A_{n,0}$.
It is now easy to establish \eqref{eq:swtilde3}.
\end{proof}

\begin{rem}\label{thm:remark}
{\rm The orthonormal polynomials $\tilde P_n(x;p,q)$ belong to a determinate moment problem. From
Theorem~\ref{thm:GSW} it is possible to find the asymptotic
behaviour of  $\tilde P_n(x;p,q)$ as $n\to\infty$ for any $x\in\mathbb C$, namely
\begin{equation}
\label{eq:asympdet}
\tilde P_n(x;p,q)\sim (-1)^nc(x)q^{-n/2},
\end{equation}
where
$$
c(x)=q^{-1/4}\frac{1-q}{1-p}\sqrt{\frac{(p;q)_\infty}{(q;q)_\infty}}\sum_{k=0}^\infty
\frac{q^{k^2+k/2}}{(pq,q;q)_k}\left(-qx\right)^k
$$
essentially is the $q$-Bessel function $J^{(2)}_{\nu}(z;q)$ with $p=q^\nu$.

To see this, we notice that
$$
\sum_{k=0}^n (-1)^k\left[\begin{matrix}n\\k\end{matrix}\right]_q
\frac{q^{k^2+{k}/{2}}}{(p;q)_k}\left[1-\frac{1-q^k}{1-pq^k}\frac{(q^{n+1};q)_\infty}{(pq^{n+1};q)_\infty}\right]x^k
$$
converges to
$$
\sum_{k=0}^\infty
(-1)^k\frac{q^{k^2+{k}/{2}}}{(p,q;q)_k}\left[1-\frac{1-q^k}{1-pq^k}\right]x^k
=\sum_{k=0}^\infty
\frac{q^{k^2+k/2}}{(pq,q;q)_k}\left(-qx\right)^k.
$$
From the $q$-binomial theorem, we find
\begin{equation}
\label{eq:asymp}
1-\frac{(q^n;q)_\infty}{(pq^n;q)_\infty} \sim \frac{1-p}{1-q}q^n
\,\mbox{ as } \, n\to\infty
\end{equation}
and combining the above, we get \eqref{eq:asympdet}.}
\end{rem}

The monic polynomials $\tilde p_n(x;p,q)=\tilde P_n(x;p,q)/\tilde b_{n,n}$
satisfy the three-term recurrence relation
\begin{equation}\label{eq:gen3t}
\tilde p_n(x;p,q)=(x-\tilde c_n)\tilde p_{n-1}(x;p,q) -
\tilde\la_n\tilde p_{n-2}(x;p,q),\quad n\ge 1,
\end{equation}
where the coefficients are given by
\begin{equation}
\label{eq:help1}
\tilde c_1=-\frac{\tilde b_{0,1}}{\tilde b_{1,1}}, \quad
\tilde c_{n+1}=\frac{\tilde b_{n-1,n}}{\tilde b_{n,n}}-\frac{\tilde
  b_{n,n+1}}{\tilde b_{n+1,n+1}}, \quad n\geq 1
\end{equation}
and
\begin{equation}
\label{eq:help2}
\tilde\la_{n+1}=\frac{\tilde b_{n-1,n-1}^2}{\tilde b_{n,n}^2},\quad
n\ge 1.
\end{equation}
Using the expressions from Theorem~\ref{thm:GSW}, we get

\begin{thm}
\label{thm:anbn}
Let $\Delta_n$ be defined as in \eqref{eq:rho}.
Then the coefficients in \eqref{eq:help1}--\eqref{eq:help2} are given by
\begin{equation}
\label{eq:cnlan}
\begin{aligned}
\tilde c_1&=\frac{(p;q)_\infty}{\Delta_1}q^{-3/2}, \\
\tilde{c}_{n+1}&=
\left[(1-q^{n+1})(pq^n;q)_\infty-(1-pq^{n+1})(q^n;q)_\infty\right]\frac{q^{-2n-3/2}}{(1-q)\Delta_{n+1}} \\
&\quad-\left[(1-q^{n})(pq^{n-1};q)_\infty-(1-pq^{n})(q^{n-1};q)_\infty\right]\frac{q^{-2n+1/2}}{(1-q)\Delta_n}
\end{aligned}
\end{equation}
and
\begin{equation}
\label{eq:lnlan}
\tilde\la_{n+1}=\frac{\Delta_{n-1}\Delta_{n+1}}{\Delta_{n}^2}(1-q^{n})(1-pq^{n})q^{-4n}.
\end{equation}
\end{thm}

\begin{proof} Specializing \eqref{eq:swtilde1} to $k=n$ and $k=n-1$, we find
$$
\tilde b_{n,n}=q^{n^2+n+1/4}\left(\frac{\Delta_n}{\Delta_{n+1}(p;q)_{n+1}(q;q)_n}\right)^{1/2}
$$
and
$$
\tilde b_{n-1,n}=
-\frac{q^{n^2-n+3/4}}{(1-q)}\frac{(1-q^n)(pq^{n-1};q)_\infty-(1-pq^n)(q^{n-1};q)_\infty}
{\sqrt{(p;q)_{n+1}(q;q)_n\Delta_n\Delta_{n+1}}}.
$$
Using \eqref{eq:help1}--\eqref{eq:help2}, we obtain the expressions in \eqref{eq:cnlan}--\eqref{eq:lnlan}.
\end{proof}

In the special case $p=q$, the formulas of Theorem~\ref{thm:GSW} and Theorem~\ref{thm:anbn}
simplify.

\begin{cor}
\label{thm:p=q}
The coefficients of \eqref{eq:swtilde} in the case $p=q$ are given by
\begin{equation}
\label{eq:swtildep=q1}
\tilde b_{k,n}=\tilde C_n (-1)^k\left[\begin{matrix}n\\k\end{matrix}\right]_q
\frac{q^{k^2+{k}/{2}}}{(q;q)_{k+1}}\left[1-q^{k+1}-(1-q^k)(1-q^{n+1})\right],
\end{equation}
where
\begin{equation}
\label{eq:swtildep=q2}
\tilde C_n=(-1)^n q^{-{n}/{2}-{1}/{4}},
\end{equation}
i.e.,
\begin{equation}
\label{eq:swtildep=q3}
\tilde b_{k,n}=b_{k,n}\,q^{-n-1/2}\left[1-\frac{(1-q^k)(1-q^{n+1})}{1-q^{k+1}}\right].
\end{equation}
Moreover,
\begin{equation}
\label{eq:swtildep=q4}
\tilde D_n=q^{n+1}D_n.
\end{equation}
Finally, the coefficients \eqref{eq:help1}--\eqref{eq:help2} in the
three-term recurrence relation are
\begin{equation}
\label{eq:3tmax}
\tilde c_n=\left(1+q^3-(1+q^2)q^n\right)q^{-2n-1/2}, \quad
\tilde \la_{n+1}=(1-q^n)^2q^{-4n}.
\end{equation}
\end{cor}

\section{The Kernel Polynomials}

By \eqref{eq:SW1}, the polynomials $S_n(x;p,q)$ are orthogonal with respect to the
density
\begin{equation}\label{eq:dens}
\mathcal D(x;p,q)=\frac{1}{\sqrt{2\pi\log(1/q)}}\exp\left(-\frac{(\log
x)^2}{2\log(1/q)}\right)(p,-p/\sqrt{q}x;q)_\infty
\end{equation}
and we see that
\begin{equation}\label{eq:densrel}
\mathcal D(qx;pq,q)=x\frac{\sqrt{q}}{1-p}\mathcal D(x;p,q).
\end{equation}
This shows that the monic polynomials $k_n(x)=q^{-n}S_n(qx;pq,q)$ are
orthogonal with respect to the density in \eqref{eq:densrel}, hence
equal to the monic kernel polynomials corresponding to $S_n(x;p,q)$.

The three-term recurrence relation for $S_n(x;p,q)$ is
$$
S_n(x;p,q)=(x-c_n)S_{n-1}(x;p,q)-\lambda_nS_{n-2}(x;p,q), \quad n\ge 1
$$
with
\begin{equation}
\label{eq:3tgenSW}
c_n=\left(1+q - (p+q)q^{n-1}\right)q^{-2n+1/2}, \quad
\lambda_{n+1}=(1-q^n)(1-pq^{n-1})q^{-4n}.
\end{equation}
It follows that the coefficients in \eqref{eq:3t} for the case
$p_n(x)=S_n(x;p,q)$ are given by
\begin{equation}
\label{eq:3tgenSW1}
d_n=\left(1+q-(1+p)q^n\right)q^{-2n-1/2},\quad
\nu_{n+1}=(1-q^n)(1-pq^n)q^{-4n-2}.
\end{equation}

Chihara observed that for $p=q$ we have the
following simple form of the coefficients in \eqref{eq:3tgenSW1}:
\begin{equation}
\label{eq:3tspec}
d_n=(1+q)(1-q^n)q^{-2n-1/2},\quad
\nu_{n+1}=(1-q^n)(1-q^{n+1})q^{-4n-2}.
\end{equation}
In this case, the chain sequence \eqref{eq:chainseq} becomes the constant
sequence
\[
\beta_n=\frac{q}{(1+q)^2}
\]
satisfying $0<\beta_n<1/4$, and the maximal parameter sequence is also constant
\[
M_n=\frac{1}{1+q}.
\]

For the shell polynomials $p_n^M$, which are equal to $\tilde
p_n(x;q,q)$, Chihara gave the following form of the coefficients from
\eqref{eq:pnh}:
$$
c_n^M=(1+q^3-(1+q^2)q^n)q^{-2n-1/2},
\quad
\la_{n+1}^M=(1-q^n)^2q^{-4n}
$$
(there is a misprint in \cite{Ch2}: The power 2 is missing in the last
formula). They agree with the expressions in \eqref{eq:3tmax}.

Going back to arbitrary $0\le p<1$, we find

\begin{thm}
\label{thm:chaingenSW}
The chain sequence \eqref{eq:chainseq} corresponding to the kernel polynomials $k_n(x)=q^{-n}S_n(qx;pq,q)$ is
\begin{equation}
\label{eq:betanSW}
\b_n=\frac{q(1-q^n)(1-pq^n)}{(1+q-(1+p)q^n)(1+q-(1+p)q^{n+1})},\quad
n\ge 1.
\end{equation}
The maximal and minimal parameter sequences $(M_n)$ and $(m_n)$ are given by
\begin{equation}\label{maxminSW}
M_n=\frac{q}{1+q-(1+p)q^{n+1}}\frac{\Delta_n}{\Delta_{n+1}},
\quad m_n=\frac{q(1-q^n)}{1+q-(1+p)q^{n+1}},
\end{equation}
and the generalized Stieltjes--Wigert polynomials $S_n(x;p,q)$ correspond to the
parameter sequence
\begin{equation}\label{eq:between}
h_n=\frac{q(1-pq^n)}{1+q-(1+p)q^{n+1}}.
\end{equation}
\end{thm}

\begin{proof} The expression for $\b_n$ follows immediately from
\eqref{eq:3tgenSW1}. We know from Theorem~\ref{thm:chain} that
$p_n^M(x)=\tilde p_n(x;p,q)$. So by \eqref{eq:pnh1},
\[
c_1^M=M_0d_1,
\]
and by \eqref{eq:cnlan} and \eqref{eq:3tgenSW1}, we have
$$
c_1^M=\tilde c_1=\frac{(p;q)_\infty}{\Delta_1}q^{-3/2},
\quad
d_1=(1+q-(1+p)q)q^{-5/2}.
$$
Hence,
$$
M_0=\frac{q(p;q)_\infty}{(1+q-(1+p)q)\Delta_1},
$$
showing the formula for $M_n$ for $n=0$. It is now easy to show by
induction that $\b_n=M_n(1-M_{n-1})$ for $n\ge 1$.

It is similarly easy to see by induction that the sequences
$(m_n),(h_n)$ are parameter sequences for $(\b_n)$. Since $m_0=0$,
it is the minimal parameter sequence.
To see that $(h_n)$ corresponds to $S_n(x;p,q)$, it suffices to verify
that $h_0d_1=c_1$, where $c_1$ is given by \eqref{eq:3tgenSW}.
\end{proof}

The parameter sequences from Theorem \ref{thm:chaingenSW} enable us to find the value $\beta$ of the continued fraction
\begin{equation}
\label{eq:cf}
1-\cfrac{\beta_1}{1-
 \cfrac{\beta_2}{1-
  \cfrac{\beta_3}{1-\dotsb
}}}
\end{equation}
in three different ways. By the results in \cite[Chap.~III]{Ch1} (see also \cite[Sect.~19]{Wa}), we have
\[
\beta=M_0=\frac{1}{1+L}=h_0+\frac{1-h_0}{1+G},
\]
where
\[
L=\sum_{n=1}^\infty\frac{m_1\cdots m_n}{(1-m_1)\cdots(1-m_n)}, \quad
G=\sum_{n=1}^\infty\frac{h_1\cdots h_n}{(1-h_1)\cdots(1-h_n)}.
\]
Since $(M_{n+k})$ is the maximal parameter sequence for the chain sequence $(\beta_{n+k})$, we can in fact find the value of
\begin{equation}
\label{eq:cf1}
1-\cfrac{\beta_{k+1}}{1-\cfrac{\beta_{k+2}}{1-\cfrac{\beta_{k+3}}{1-\dotsb}}}
\end{equation}
for every $k\geq 0$.

We collect the above considerations in
\begin{thm}
Let $(\beta_n)$ be the chain sequence given by \eqref{eq:betanSW}. Then the continued fraction in \eqref{eq:cf}
has the value
\[
\beta=\frac{q}{1-pq}\frac{\Delta_0}{\Delta_1}=\frac{q(1-p)(pq^2;q)_\infty}{(pq;q)_\infty-(q;q)_\infty}.
\]
More generally, the continued fraction in \eqref{eq:cf1} has the value
\[
M_k=\frac{q}{1+q-(1+p)q^{k+1}}\frac{\Delta_{k}}{\Delta_{k+1}}, \quad k\geq 0.
\]
\end{thm}
\begin{proof}
The result follows immediately from \cite[Thm.~6.1 (Chap.~III)]{Ch1}. To find $L$ and $G$, note that
\[
\frac{m_k}{1-m_k}=\frac{q(1-q^k)}{1-pq^{k+1}},
\quad
\frac{h_k}{1-h_k}=\frac{q(1-pq^k)}{1-q^{k+1}},
\]
so that
\[
1+L=\sum_{n=0}^\infty\frac{(q;q)_n}{(pq^2;q)_n}q^n,
\quad
1+G=\sum_{n=0}^\infty\frac{(pq;q)_n}{(q^2;q)_n}q^n.
\]
The value of $1+G$ can thus be found using the $q$-binomial theorem.
To compute $1+L$, one first applies Heine's transformation formula and then the $q$-binomial theorem.
\end{proof}

\begin{rem}
{\rm We mention that
\begin{equation}
\label{eq:Max}
\sum_{n=1}^\infty\frac{M_1\cdots M_n}{(1-M_1)\cdots(1-M_n)}=\infty,
\end{equation}
precisely as should be the case for the maximal parameter sequence. To see this, note that
\[
\frac{M_k}{1-M_k}=\frac{M_k M_{k+1}}{\beta_{k+1}}=
\frac{\Delta_k}{\Delta_{k+2}}\frac{q}{(1-q^{k+1})(1-pq^{k+1})}
\]
so that the series in \eqref{eq:Max} reduces to
\[
\sum_{n=1}^\infty\frac{\Delta_1\Delta_2}{\Delta_{n+1}\Delta_{n+2}}\frac{q^n}{(q^2,pq^2;q)_n}.
\]
On the lines of \eqref{eq:asymp}, we have
\[
\Delta_n=\frac{1-p}{1-q}q^n+\mathcal{O}(q^{2n})
\]  
and the result follows.}
\end{rem}

\noindent
Christian Berg, Jacob Stordal Christiansen\\
Department of Mathematics, University of Copenhagen,
Universitetsparken 5, DK-2100, Denmark\\
e-mail: {\tt berg@math.ku.dk}\\
e-mail: {\tt stordal@math.ku.dk}


\vspace{0.8cm}
\begin{thebibliography}{120}

\bibitem{Ak} N.~I.~Akhiezer, {\it The Classical Moment Problem and Some Related
 Questions in Analysis}. English translation, Oliver and Boyd, Edinburgh,
 1965.

\bibitem{B:C}  C.~Berg, J.~P.~R.~Christensen, {\it Density questions
    in the classical theory of moments}, Ann. Inst. Fourier, Grenoble
  {\bf 31},3 (1981), 99--114.

\bibitem{B:D}  C.~Berg, A.~J.~Dur{\'a}n, {\it The index of
determinacy for measures and the $\ell ^2$-norm of orthogonal polynomials},
Trans. Amer. Math. Soc. {\bf 347} (1995), 2795--2811.

\bibitem{B:S} C.~Berg, R.~Szwarc, {\it The smallest Eigenvalue of
    Hankel Matrices}, Constr. Approx. (To appear) DOI 10.1007/s00365-010-9109-4.

\bibitem{B:T} C.~Berg, M.~Thill, {\it Rotation invariant moment
    problems}, Acta Math. {\bf 167} (1991), 207--227.

\bibitem{B:V} C.~Berg, G.~Valent, {\it The Nevanlinna parametrization
for some indeterminate Stieltjes moment problems associated with birth and
death processes}, Methods and Applications of Analysis {\bf 1} (1994),
169--209.

\bibitem{B:W} W.~Bergweiler, W.~K.~Hayman,
{\it Zeros of solutions of a functional equation,}
Comput. Methods Funct. Theory {\bf 3} (2003), 55--78.

\bibitem{Ch} T.~S.~Chihara, {\it Chain sequences and orthogonal polynomials,} Trans.
Amer. Math. Soc. {\bf 104} (1962), 1--16.

\bibitem{Ch1} T.~S.~Chihara, {\it An Introduction to Orthogonal
    Polynomials}. Gordon and Breach, New York 1978.

\bibitem{Ch2} T.~S.~Chihara, {\it Shell polynomials and indeterminate
    moment problems}, J. Comput. Appl. Math. {\bf 133} (2001), 680--681.

\bibitem{Chr} J.~S.~Christiansen, {\it The moment problem associated
    with the Stieltjes--Wigert polynomials}, J. Math. Anal. Appl. {\bf 277} (2003),
218--245.

\bibitem{G:R} G.~Gasper, M.~Rahman, {\it Basic hypergeometric series}.
Cambridge University Press, Cambridge 1990, second edition 2004.

\bibitem{Hay} W.~K.~Hayman, {\it On the zeros of a {$q$}-{B}essel function},
 Contemp. Math. {\bf 382} (2005), 205--216.

\bibitem{Hub} T.~Huber, {\it Hadamard products for generalized {R}ogers-{R}amanujan series,}
J. Approx. Theory {\bf 151} (2008), 126--154.

\bibitem{I:Z} M.~E.~H.~Ismail, C.~Zhang,
{\it Zeros of entire functions and a problem of {R}amanujan,}
{Adv. Math.} {\bf 209} (2007), 363--380.

\bibitem{S:T} J.~Shohat, J.~D.~Tamarkin, {\it The Problem of Moments.}
Revised edition, American Mathematical Society, Providence, 1950.

\bibitem{St} T.J.\ Stieltjes, {\it Recherches sur les fractions continues},
Annales de la Facult\'e des Sciences de Toulouse, {\bf 8} (1894),
1--122; {\bf 9} (1895), 5--47. English translation in Thomas Jan
Stieltjes, {\it Collected papers}, Vol. II,
pp. 609--745. Springer-Verlag, Berlin, Heidelberg. New York, 1993.

\bibitem{Sz1} G.~Szeg{\H o}, Orthogonal Polynomials, 4th ed.,
  Colloquium  Publications, vol. 23, Amer. Math. Soc., Rhode Island, 1975.

\bibitem{Wa} H.~S.~Wall, Analytic {T}heory of {C}ontinued {F}ractions,
 D.~Van Nostrand Company, Inc., New York, N. Y., {1948}.

\bibitem{Wig} S.~Wigert, {\it Sur les polyn{\^o}mes orthogonaux et
    l'approximation des fonctions continues}, Arkiv f\"or Matematik,
  Astronomi och Fysik {\bf 17} (1923), no. 18, 15pp.

\end{thebibliography}
\end{document}